\documentclass[12pt]{article}
\def\be{\begin{equation}}
\def\ee{\end{equation}}
\def\ff#1{\mbox{\boldmath $#1$} }
\def\a{\alpha}
\def\b{\beta}
\def\lam{\lambda}
\def\e{\epsilon}
\def\x{\ff{x}}

\def\Q{{\Omega}}
\def\bP{\ff{P}}
\newcommand{\mat}[1]{{\left( \begin{array}{cccc}#1\end{array}\!\right)}}

\begin{document}

\title{\Large \bf Free Lunch or No Free Lunch: \\ That is not Just a Question?}

\author{Xin-She Yang \\
Mathematics and Scientific Computing \\ National Physical Laboratory, Teddington, TW11 0LW, UK}

\date{}

\maketitle

\begin{abstract}
The increasing popularity of metaheuristic algorithms has attracted a great deal of attention
in algorithm analysis and performance evaluations. No-free-lunch theorems are of both theoretical
and practical importance, while many important studies on
convergence analysis of various metaheuristic algorithms have proven to be fruitful.
This paper discusses the recent results on no-free-lunch theorems and algorithm convergence,
as well as their important implications for algorithm development in practice.
Free lunches may exist for certain types of problem. In addition, we will highlight
some open problems for further research. \\

\noindent {\bf Citation Details:} X. S. Yang, Free lunch or no free lunch: that is not just a question?
{\it Int. J. Artificial Intelligence Tools}, Vol. 21, No. 3,  1240010 (2012).  \\
DOI: 10.1142/S0218213012400106 

\end{abstract}

\newpage

\section{Introduction}

Metaheuristic algorithms form an important part of contemporary global
optimization algorithms \cite{Auger,Auger2,Blum,Yang,Neumann,Parp}. Good examples are simulated annealing
and particle swarm optimization \cite{Kennedy,Kirk}.
They work remarkably efficiently and have many
advantages over traditional, deterministic methods and algorithms, and
thus they have been applied in almost all area of science, engineering
and industry \cite{Cui,Floudas,Cui2,Yang4}.

Despite such a huge success in applications, mathematical analysis of
algorithms remains limited and many open problems are still un-resolved.
There are three challenging areas for algorithm analysis: complexity,
convergence and no-free-lunch theory.

Complexity analysis of traditional algorithms such as quick sort and
matrix inverse are well-established, as these algorithms are deterministic.
In contrast, complexity analysis of metaheuristics remains a challenging
task, as the stochastic nature of these algorithms. However, good results
do exist, concerning randomization search techniques \cite{Auger}.

Convergence analysis is another challenging area. One of the main
difficulties concerning the convergence analysis of
metaheuristic algorithms is the lack of a generic framework, though
substantial studies have been carried out using dynamic systems and
Markov processes. However, convergence analysis still remains
one of the active research areas with many encouraging results
\cite{Clerc,Trelea,Ola,Gutja}.

Along the mathematical analysis of optimization algorithms,
another equally challenging, and yet fruitful area is the
theory on algorithm performance and comparison, leading
to a wide range of no-free-lunch (NFL)
theorems \cite{Wolpert,Igel}.
While in well-posed cases of optimization functional space in finite domains, NFL theorems
do hold, however, free lunches are possible \cite{Auger,Wolpert2,Villa}.

In this paper, we will briefly review and summarize the recent
studies of no-free-lunch theory and also free lunch scenario.
This enable us to view the NLF and free lunch in a unified
framework, or at least, in a convenient way. We will also
briefly highlights some of the convergence studies.
Based on these studies, we will summarize and propose
a series of recommendations for further research.

\section{No-Free-Lunch Theorems}

The seminal paper by Wolpert and Mcready in 1997 essentially
proposed a framework for performance comparison of
optimization algorithms, using a combination of Bayesian statistics
and Markov random field theories.

Along many relevant assumptions in proving the NFL theorems, two
fundamental assumptions are: finite states of the search space
(and thus the objective values), and the non-revisiting time-ordered sets.

The first assumption is a good approximation to many problems, especially
in finite-digit approximations. However, there is mathematical difference
in countable finite, and countable infinite. Therefore, the results
for finite states/domains may not directly applicable to infinite domains.
Furthermore, as continuous problem are uncountable, NFL results for finite
domains will usually not hold for continuous domains \cite{Auger}.

The second assumption on non-revisiting iterative sequence is an over-simplification,
as almost all metaheuristic algorithms are revisiting in practice, some points
visited before will possibly be re-visited again in the future. The only possible
exception is the Tabu algorithm with a very long Tabu list \cite{Glover}.
Therefore, results for non-revisiting time-ordered iterations may not be true
for the cases of revisiting cases, because the revisiting iterations break an
important assumption of `closed under permutation' (c.u.p) required for proving
the NFL theorems \cite{Marshall}.

Furthermore, optimization problems do not necessarily concern the whole set of all possible
functions/problems, and it is often sufficient to consider a subset of problems.
It is worth pointing out active studies have carried out in constructing algorithms that
can work best on specific subsets of optimization problems, in fact, NFL theorems do not
hold in this case \cite{Christensen}.

Before we go further to discuss more about any possible free lunches, let us
sketch Wolpert and Macready's original proof. Assuming that the
search space is finite (though quite large), thus the space of possible objective
values is also finite. This means that objective function is simply a mapping
$ f: {\cal X} \mapsto {\cal Y}$, with  ${\cal F}={\cal Y}^{\cal X}$ as the
space of all possible problems under permutation.

As an algorithm tends to produce a series of points or solutions
in the search space, it is further assumed that these points are distinct.
That is,  for $k$ iterations, $k$ distinct visited points forms a
time-ordered set
\be \Q_k=\Big\{\Big(\Q_k^x(1),\Q_k^y(1)\Big),...,\Big(\Q_k^x(k), \Q_k^y(k)\Big)\Big\}. \ee

There are many ways to define a performance measure, though a good measure
still remains debatable \cite{Shilane}. Such a measure can depend on
the number of iteration $k$, the algorithm $a$ and the actual cost function $f$, which
can be denoted by $ P(\Q_k^y\|f,k,a)$. Here we follow the notation style in seminal paper
by Wolpert and Mcready (1997). For any pair of algorithms $a$ and $b$, the NFL theorem states
\be \sum_f P(\Q_k^y|f,k,a) =\sum_f P(\Q_k^y |f,k,b). \ee
In other words, any algorithm is as good (bad) as a random search, when
the performance is averaged {\emph over} all possible functions.

Wolpert and Macready's original proof was carried out by induction. Using a similar methodology
and similar assumptions, other forms of NFL theorems
can also be proved. These theorems are vigorous and thus have important theoretical
values. However, their practical implications are a different issue. In fact,
it may not be so important in practice anyway, we will discuss this in a later section.

\section{Free Lunch or No Free Lunch}

\subsection{Continuous Free Lunches}

The validity of NFL theorems largely depends on the validity of their fundamental
assumptions. However, whether these assumptions are valid in practice is
another question. Often, these assumptions are too stringent, and thus free lunches are possible.

One of the assumptions is the non-revisiting nature of the $k$ distinct points
which form a time-ordered set. For revisiting points as they do occur in
practice in real-world optimization algorithms,  the `closed under permutation'
does not hold, which renders NFL theorems invalid  \cite{Schumacher,Marshall}.
This means free lunches do exist in practical applications.

Another basic assumption is the finiteness of the domains. For continuous domains,
Auger and Teytaud in 2010 have proven that the NFL theorem does not hold, and therefore
they concluded that ``continuous free lunches exist". Indeed,
some algorithms are better than others. For example, for a 2D sphere function,
they demonstrated that an efficient algorithm only needs 4 iterations/steps
to reach the global minimum.

\subsection{Coevolutionary and Multiobjective Free Lunches}

The basic NFL theorems concern a single agent, marching iteratively in the search
space in distinct steps. However, Wolpert and Mcready proved in 2005 that
NFL theorems do not hold under coevolution.  For example, a set of players (or agents)
in self-play problems can work together so as to produce a champion.
This can be visualized as an evolutionary process of training a chess champion.
In this case, free lunch does exist \cite{Wolpert2}.
It is worth pointing out that for a single player, it tries to pursue the
best next move, while for two players, the fitness function depend on the moves
of both players. Therefore, the basic assumptions for NFL theorems are no longer valid.

For multiobjective optimization problems, things have become even more complicated.
An important step in theoretical analysis is that some multiobjective optimizers are better than others
as pointed out by Corne and Knowles \cite{Corne}. One of the major reasons is that
the archiver and generator in the multiobjective algorithms can lead to
multiobjective free lunches.

Whether NFL holds or not, it has nothing to say about the complexity of the problems.
In fact,  no free lunch theorem has not been proved to be true for problems
with NP-hard complexity \cite{Whitley}.

\section{Practical Implications of NFL Theorems}


No-free-lunch theorems may be of theoretical importance, and they can also
have important implications for algorithm development in practice, though
not everyone agrees the real importance of these implications.

There are three kinds of opinions concerning the implications. The first group
may simply ignore these theorems, as they argue that the assumptions are
too stringent, and the performance measures based on average overall functions are irrelevant
in practice. Therefore, no practical importance can be inferred, and research just carries on.

The second kind is that NFL theorems can be true, and they can accept that the fact
there is no universally efficient algorithm. But in practice some algorithms do
performance better than others for a specific problem or a subset of problems.
Research effort should focus on finding the right algorithms for the right type of problem.
Problem-specific knowledge is always helpful to find the right algorithm(s).

The third kind of opinion is that NFL theorems are not true for other types of problems such
as continuous problems and NP-hard problems. Theoretical work concerns more elaborate
studies on extending NFL theorems to other cases or on finding free lunches \cite{Auger}.
On the other hand, algorithm development continues to design better algorithms
which can work for a wider range of problems, not necessarily all types of problems.
As we have seen from the above analysis, free lunches do exist, and better algorithms
can be designed for a specific subset of problems \cite{Yang2,YangDeb}.

Thus, free lunch or no free lunch is not just a simple question, it has important and yet practical
importance. There is certain truth in all the above arguments, and their impacts on optimization community
are somehow mixed. Obviously, in reality, the algorithms with problem-specific knowledge
typically work better than random search, and we also realized that there is no
universally generic tool that works best for all the problems. Therefore,
we have to seek balance between
speciality and generality, between algorithm simplicity and problem complexity, and between problem-specific
knowledge and capability of handling black-box optimization problems.

\section{Convergence Analysis}

For convergence analysis, there is no mathematical framework in general to provide
insights into the working mechanism, the complexity, stability and convergence
of any given algorithm \cite{He,Thiko}.
Despite the increasing popularity of metaheuristics, mathematical analysis
remains fragmental, and many open problems concerning convergence analysis
need urgent attention.

\subsection{PSO}

Particle swarm optimization (PSO) was developed by Kennedy and
Eberhart in 1995 \cite{Kennedy}, based on the swarm behaviour such
as fish and bird schooling in nature. Since then, PSO has
generated much wider interests, and forms an exciting, ever-expanding
research subject, called swarm intelligence. PSO has been applied
to almost every area in optimization, computational intelligence,
and design/scheduling applications.

The movement of a swarming particle consists of two major components:
a stochastic component and a deterministic component.
Each  particle is attracted toward the position of the current global best
$\ff{g}^*$ and its own best location $\x_i^*$ in history,
while at the same time it has a tendency to move randomly.

Let $\x_i$ and $\ff{v}_i$ be the position vector and velocity for
particle $i$, respectively. The new velocity and location updating formulas are
determined by
\be \ff{v}_i^{t+1}= \ff{v}_i^t  + \a \ff{\e}_1
[\ff{g}^*-\x_i^t] + \b \ff{\e}_2 [\x_i^*-\x_i^t].
\label{pso-speed-100}
\ee
\be \x_i^{t+1}=\x_i^t + \ff{v}_i^{t+1}, \label{pso-speed-140} \ee
where $\ff{\e}_1$ and $\ff{\e}_2$ are two random vectors, and each
entry taking the values between 0 and 1. The parameters $\a$ and $\b$ are the learning parameters or
acceleration constants, which can typically be taken as, say, $\a \approx \b \approx 2$.

There are at least two dozen PSO variants which extend the standard PSO
algorithm, and the most noticeable improvement
is probably to use inertia function $\theta
(t)$ so that $\ff{v}_i^t$ is replaced by $\theta(t) \ff{v}_i^t$
where $\theta \in [0,1]$.
This is equivalent to introducing a virtual mass to stabilize the motion
of the particles, and thus the algorithm is expected to converge more quickly.

The first convergence analysis of PSO was carried out by Clerc and Kennedy
in 2002 \cite{Clerc} using the theory of dynamical systems.
Mathematically, if we ignore the random factors, we can view the system
formed by (\ref{pso-speed-100}) and (\ref{pso-speed-140}) as a dynamical system.
If we focus on a single particle $i$ and imagine that there is only one particle in this system,
then the global best $\ff{g}^*$ is the same as its current best $\x_i^*$. In this case, we have
\be \ff{v}_i^{t+1} = \ff{v}_i^t + \gamma (\ff{g}^*-\x_i^t), \quad \gamma=\a+\b, \ee
and
\be \x_i^{t+1} =\x_i^t + \ff{v}_i^{t+1}. \ee
Considering the 1D dynamical system for particle swarm optimization,
we can replace $\ff{g}^*$ by a parameter constant $p$ so
that we can see if or not the particle of interest will converge towards $p$.
By setting $u_t=p-x(t+1)$
and using the notations for dynamical systems, we have a simple dynamical system
\be v_{t+1}=v_t + \gamma u_t, \quad u_{t+1} =-v_t + (1-\gamma) u_t, \ee
or
\be Y_{t+1} = A Y_t, \quad A=\mat{ 1 & & \gamma \\ -1 & & 1-\gamma}, \quad Y_t=\mat{v_t \\ u_t}. \ee
The general solution of this dynamical system can be written as $Y_t=Y_0 \exp[A t]$.
The system behaviour can be characterized by the eigenvalues $\lam$ of $A$
\be \lam_{1,2}=1-\frac{\gamma}{2} \pm \frac{\sqrt{\gamma^2 - 4 \gamma}}{2}. \ee
It can be seen clearly that $\gamma=4$ leads to a  bifurcation.

Following a straightforward analysis of this dynamical system,
we can have three cases. For $0 < \gamma <4$, cyclic and/or quasi-cyclic
trajectories exist. In this case, when randomness is gradually reduced,
some convergence can be observed.
For $\gamma>4$, non-cyclic behaviour can be expected and
the distance from $Y_t$ to the center $(0,0)$ is monotonically increasing with $t$.
In a special case $\gamma=4$, some convergence behaviour can be observed.
For detailed analysis, please refer to Clerc and Kennedy \cite{Clerc}.
Since $p$ is linked with the global best, as the iterations continue,
it can be expected that all particles will aggregate towards the
the global best.

\subsection{Firefly Algorithm}

Firefly Algorithm (FA) was developed by Yang \cite{Yang,Yang3},
which was based on the flashing patterns and behaviour
of fireflies. In essence, each firefly will be attracted to brighter ones,
while at the same time, it explores and searches for prey randomly.
In addition, the brightness of a firefly is determined by the landscape of the
objective function.

The movement of a firefly $i$ is attracted to another more attractive (brighter)
firefly $j$ is determined by
\be \x_i =\x_i + \b_0 e^{-\gamma r^2_{ij}} (\x_j-\x_i) + \a \; \ff{\e}_i,
\label{FA-equ-50}  \ee
where the second term is due to the attraction. The third term
is randomization with $\a$ being the randomization parameter, and
$\ff{\e}_i$ is a vector of random numbers drawn from a Gaussian distribution
or uniform distribution.  Here is $\b_0 \in [0,1]$ is the attractiveness at $r=0$,
and $r_{ij}=||\x_i-\x_j||$ is the Cartesian distance. For other problems such as scheduling,
any measure that can effectively characterize the quantities of interest in the optimization
problem can be used as the `distance' $r$.
For most implementations, we can take $\b_0=1$, $\a=O(1)$ and $\gamma=O(1)$.
It is worth pointing out that (\ref{FA-equ-50}) is essentially a random walk biased
towards the brighter fireflies. If $\b_0=0$, it becomes a simple random walk.
Furthermore, the randomization term can easily be
extended to other distributions such as L\'evy flights.

We now can carry out the convergence analysis for the firefly algorithm in a
framework similar to Clerc and Kennedy's dynamical analysis. For simplicity,
we start from the equation for firefly motion without the randomness term
\be \x_i^{t+1}=\x_i^t + \b_0 e^{-\gamma r_{ij}^2} (\x_j^t - \x_i^t). \ee
If we focus on a single agent, we can replace $\x_j^t$ by the global best
$g$ found so far, and we have
\be \x_i^{t+1} =\x_i^t + \b_0 e^{-\gamma r_{i}^2 } (g-\x_i^t), \ee
where the distance $r_i$ can be given by the $\ell_2$-norm $r_i^2=||g-\x_i^t||_2^2$.
In an even simpler 1-D case, we can set $y_{t}=g-\x_i^t$, and we have
\be y_{t+1} =y_t - \b_0 e^{-\gamma y_t^2} y_t. \label{fa-dynamics-25} \ee
We can see that $\gamma$ is a scaling parameter which only affects the scales/size
of the firefly movement. In fact, we can let
$u_t=\sqrt{\gamma} y_t$ and we have
\be u_{t+1}=u_t [1-\b_0 e^{-u_t^2}]. \label{fa-dynamics-150} \ee
These equations can be analyzed easily using the same methodology for studying
the well-known logistic map \be u_t=\lam u_t (1-u_t). \label{fa-chaos-55} \ee

Straightforward analysis can show that convergence can be achieved
for $\b_0 <2$. There is a transition from periodic to chaos at $\b_0 \approx 4$.
This may be surprising, as the aim of designing a metaheuristic algorithm is
to try to find the optimal solution efficiently and accurately. However,
chaotic behaviour is not necessarily a nuisance, in fact, we can use it to
the advantage of the firefly algorithm. Simple chaotic characteristics from  (\ref{fa-chaos-55})
can often be used as an efficient mixing technique for generating diverse solutions.
Statistically, the logistic mapping (\ref{fa-chaos-55}) for $\lam=4$ for the initial
states in (0,1) corresponds a beta distribution
\be B(u,p,q)=\frac{\Gamma(p+q)}{\Gamma(p) \Gamma(q)} u^{p-1} (1-u)^{q-1}, \ee
when $p=q=1/2$. Here $\Gamma(z)$ is the Gamma function
\be \Gamma(z) =\int_0^{\infty} t^{z-1} e^{-t} dt.  \ee
In the case when  $z=n$ is an integer, we have $\Gamma(n)=(n-1)!$. In addition,
$\Gamma(1/2)=\sqrt{\pi}$. From the algorithm implementation point of view,
we can use higher attractiveness $\b_0$ during the early stage
of iterations so that the fireflies can explore, even chaotically,
the search space more effectively. As the search continues and convergence
approaches, we can reduce the attractiveness $\b_0$ gradually, which may
increase the overall efficiency of the algorithm. Obviously, more
studies are highly needed to confirm this.

\subsection{Markov Chains}

Most theoretical studies use Markov chains/process as a framework for convergence analysis.
A Markov chain is said be to regular if some positive power $k$  of the transition matrix $\bP$
has only positive elements. A chain is ergodic or irreducible if it is aperiodic and positive recurrent,
which means that it is possible to reach every state from any state.

For  a time-homogeneous chain as $k \rightarrow \infty$, we have the stationary probability distribution $\pi$,
satisfying \be \pi=\pi \bP, \ee
thus the first eigenvalue is always $1$. This will lead to the
 asymptotic convergence to the global optimality $\theta_*$:
 \be \lim_{k \rightarrow \infty} \theta_k \rightarrow \theta_*, \ee
 with probability one \cite{Gamer,Ola,Gutja}.

Now if look at the PSO closely using the framework of Markov chain
Monte Carlo, each particle in PSO  essentially
forms a Markov chain, though this Markov chain is biased towards to the
current best, as the transition probability often leads to the
acceptance of the move towards the current global best.
Other population-based algorithms can also be viewed in this framework.
In essence, all metaheuristic algorithms with piecewise, interacting paths
can be analyzed in the general framework of Markov chain Monte Carlo.
The main challenge is to realize this and to use the appropriate Markov chain
theory to study metaheuristic algorithms. More fruitful studies will surely emerge in the future.

\subsection{Convergence of SA}

Simulated annealing and generalize hill-climber algorithms were among the first
algorithms with important results on convergence analysis
\cite{Mitra,Granville,Jacobson,YangL,Stein}.
The main idea is to consider simulated annealing as sequence of homogeneous Markov
chains or a long single inhomogeneous Markov chain \cite{Hend}.
Under weak ergodic conditions, the temperature $T_k$ can be reduced sufficiently slow to
zero by
\be T_k \ge \frac{A}{\log (k)}, \quad \lim_{k \rightarrow \infty} T_k \rightarrow 0, \ee
where $A$ is a constant

\subsection{Convergence of GA}

One of the well-studied and most popular algorithms is the class of
genetic algorithms \cite{Holland,Neumann}.
Earlier seminal papers proved the convergence of
genetic algorithms \cite{Hartl,Rudolph}.
Important studies on convergence analysis of GA have been
carried out by Aytug et al. \cite{Aytug,Aytug2},
Greenhalgh and Marshal \cite{Green}, Gutjahr \cite{Gutj,Gutja}and others.

For example, a well-known result is that
the number of iterations $t(\zeta)$ in GA with a convergence probability of $\zeta$
can be estimated by the upper limit
$$ t(\zeta) \le \Bigg\lceil \frac{\ln(1-\zeta)}{\ln \bigg\{1-\min[(1-\mu)^{L n}, \mu^{Ln}] \bigg\}}  \Bigg\rceil,$$
where parameter $\mu$ is the mutation rate in genetic algorithms.
$L$ and $n$ are the string length and population size, respectively \cite{Aytug}.
These results are further elaborated by others \cite{Green}.

\subsection{Multiobjective Metaheuristics}

Convergence analysis for single-objective optimization tends to be challenging, this complexity
is further complicated by the Pareto optimality of multiobjective optimization. Despite
these challenges, asymptotic convergence of metaheuristic, multiobjective optimization
has been proved by Villalobos-Arias et al. (2005) using a framework of Markov chains \cite{Villa}.
They proved that the transition matrix $\bP$ of a metaheuristic algorithm can have
a stationary distribution $\pi$ such that
$$ |P_{ij}^k -\pi_j| \le (1-\zeta)^{k-1}, \quad \forall i,j, \quad (k=1,2,...),$$
where $\zeta$ is a function of mutation probability $\mu$, string length $L$ and population size $n$.
For example, when the population can be divided into two sets
with mutation rates $\mu_1, \mu_2$ and population sizes $n$, $n_1$, respectively,
this $\zeta$ function becomes
\be \zeta=2^{nL} \mu_1^{n_1 L} \mu_2^{(n-n_1)L}. \ee
They demonstrated that an algorithm satisfying this condition may not converge for
multiobjective optimization problems, however, an algorithm with {\it elitism}
indeed converges under the above conditions.

\subsection{Other results}

Limited results on convergence analysis exist, concerning finite domains,
 ant colony optimization \cite{Gutj,Seba},
 cross-entropy optimization, best-so-far convergence \cite{Margolin,Gutja},
 nested partition method, Tabu search, and
 largely combinatorial optimization.
However, more challenging tasks for infinite states/domains and continuous problems.
Many, many open problems need satisfactory answers.

On the other hand, it is worth pointing out that an algorithm can converge,
but it may not be efficient, as its convergence rate could be typically low.
One of the main tasks in research is to find efficient algorithms for
a given type of problem.

\section{Open Problems}

Active research on NFL theorems and algorithm convergence analysis has led to
many important results. Despite this, many crucial problems remain unanswered.
These open questions span a diverse range of areas. Here we highlight a few but
relevant open problems.

{\it Framework:} Convergence analysis has been fruitful, however, it is still
highly needed to  develop a unified framework for algorithmic analysis and convergence.

{\it Exploration and exploitation:} Two important components of metaheuristics
are exploration and exploitation or diversification and intensification.
What is the optimal balance between these two components?

{\it Performance measure:} To compare two algorithms, we have to define a measure
for gauging their performance \cite{Spall}.
At present, there is no agreed performance measure,
but what are the best performance measures ? Statistically?

{\it Free lunches:} No-free-lunch theorems have not been proved
for continuous domains for multiobjective optimization. For single-objective optimization,
free lunches are possible, is this true for multiobjective optimization?
In addition,  no free lunch theorem has not been proved to be true for problems
with NP-hard complexity (Whitley and Watson 2005).
If free lunches exist, what are their implications in practice and
how to find the best algorithm(s)?

{\it Knowledge:} Problem-specific knowledge always helps to find an appropriate solution?
How to quantify such knowledge?

{\it Intelligent algorithms:} A major aim for algorithm development is to design better,
intelligent algorithms for solving tough NP-hard optimization problems. What do mean
by `intelligent'? What are the practical ways  to design truly intelligent, self-evolving
algorithms?

\end{document}